# Multi-objective Active Control Policy Design for Commensurate and Incommensurate Fractional Order Chaotic Financial Systems

Indranil Pan[a], Saptarshi Das[b] and Shantanu Das[c]

a) Energy, Environment, Modelling and Minerals ($E^2M^2$) Research Section, Department of Earth Science and Engineering, Imperial College London, Exhibition Road, London SW7 2AZ, United Kingdom.

b) Communications, Signal Processing and Control (CSPC) Group, School of Electronics and Computer Science, University of Southampton, Southampton SO17 1BJ, United Kingdom.

c) Reactor Control Division, Bhabha Atomic Research Centre, Mumbai-400085, India.

**Authors' Emails:**

indranil.jj@student.iitd.ac.in, i.pan11@imperial.ac.uk (I. Pan)

saptarshi@pe.jusl.ac.in, s.das@soton.ac.uk (S. Das)

shantanu@magnum.barc.gov.in (Sh. Das)

**Abstract:**

In this paper, an active control policy design for a fractional order (FO) financial system is attempted, considering multiple conflicting objectives. An active control template as a nonlinear state feedback mechanism is developed and the controller gains are chosen within a multi-objective optimization (MOO) framework to satisfy the conditions of asymptotic stability, derived analytically. The MOO gives a set of solutions on the Pareto optimal front for the multiple conflicting objectives that are considered. It is shown that there is a trade-off between the multiple design objectives and a better performance in one objective can only be obtained at the cost of performance deterioration in the other objectives. The multi-objective controller design has been compared using three different MOO techniques *viz*. Non Dominated Sorting Genetic Algorithm-II (NSGA-II), epsilon variable Multi-Objective Genetic Algorithm (ev-MOGA), and Multi Objective Evolutionary Algorithm with Decomposition (MOEA/D). The robustness of the same control policy designed with the nominal system settings have been investigated also for gradual decrease in the commensurate and incommensurate fractional orders of the financial system.

## 1. Introduction

Investigations into chaotic dynamics of physical systems have revealed a variety of different fields where this is found and financial systems have been documented to show significant chaotic behaviour [1]. On contrary, fractional calculus driven modelling techniques especially fractional Brownian motion have received huge focus as a potential tool to describe the dynamical behaviour of the stochastic variations in financial time series [2]. Data driven modelling of financial systems has been shown to obey a power law characteristics i.e. the Fourier transform spectra decays as a power law with respect to frequency [3], [4]. It has been shown by Meerschaert and Scalas [5] that in finance, the relation between random variables like log-returns and waiting time can be suitably modelled using FO partial differential equations. FO noise characteristics have also been used to identify the economic periods of crisis from financial time series in [6]. Effect of parameter switching on such FO chaotic financial systems have been studied by Danca *et al.* [7]. Other perspectives of financial modelling e.g. FO volatility model has also been developed [8] for empirical market data. The concept of FO financial model has been extended to variable order financial systems [9] where the fractional orders changes over time. These studies show that sudden big fluctuations in financial time series give rise to the power law characteristics and has a close relation with fractional calculus. A realistic FO macroeconomic model was estimated using the national economic data of UK, Canada and Australia in the studies by Skovranek *et al.* [10]. Similar nonlinear model parameter estimation has been proposed using least squares to model macroeconomic data of USA [11] and interest rate change in Japan [12]. Studies have also found that the presence of time delays in such financial systems modifies the chaotic behaviour of the system where one policy change take some time to modify the overall system's dynamics [13].

It has also been found that complex financial systems show both stochastic and deterministic dynamics where the first branch has emerged to model typical behaviours like non-stationarity, non-Gaussianity, randomness and long range dependence (or power law characteristics) of such systems as discussed earlier. The other branch has emerged while analysing significant nonlinear dynamical behaviour like chaos, bifurcation [14] and hyper-chaos [15] in such large scale financial systems. There have been attempts to investigate chaotic dynamics in financial time series using delay embedding based phase space reconstruction, Lyapunov exponent estimation by parametric and non-parametric methods [16], recurrence plots [17] etc. Apart from the practical data or time series based studies, continuous time [14] as well as discrete time models [18] have been proposed to model chaotic dynamics of financial systems. Thus the co-existence of chaotic and FO characteristics are inherent in financial systems which motivates the study of an active control policy for such systems.

These chaotic dynamics are undesirable and must be supressed to reduce financial risks and improve the performance of the economy [17]. Classically there exists two broad methods for chaos control, *viz.* the OGY (Ott-Greborgi-Yorke) method of intermittent control and the continuous control method [19], [20]. FO economical or financial system [21] has been controlled or synchronized using several approaches e.g. sliding mode [22], time delayed feedback [23], linear control [24], lag projective synchronization [25], Lyapunov linearization and stability condition [26] etc. In [27], the control of the uncertain FO financial system has been attempted using adaptive sliding mode control. However in all the above cases, chaos control has been done from a stability point of view but the control performance has not been taken into consideration. Other computational intelligence based techniques which use intelligent algorithms for chaos control [28], [29] or synchronisation [30] take the performance measure of fast synchronisation or control in the formulation of the objective function itself. However, the drawback of this type of design methodology is that guaranteed analytical stability is not enforced in the process and thus the scheme may not work for initial conditions other than the one used in the simulation. Also only a single objective has been considered as a performance measure in the designs reported in [29], [30]. In a practical design problem there exists multiple trade-offs among a set of conflicting objectives. Therefore a design methodology must take these challenges into account and come up with optimal solutions which meet these objectives to a sufficient level. In other words there is a requirement for multi-objective optimisation methods to be applied to these problems to arrive at efficient designs.

Multi objective synchronisation for chaotic systems has been recently investigated in [31] where the coupling strengths between the two chaotic systems are optimised using an evolutionary multi-objective optimisation. However, in this case, the analytical stability is not included within the optimisation algorithm. Thus using the methodology proposed in [31], it might happen that for different values of initial conditions, the chaotic systems do not synchronise. This is because, the synchronisation has only been achieved in a mean squared sense and guaranteed analytical stability of the error dynamical system is not enforced. In the present paper, the concept of multi-objective synchronisation of chaotic systems in [31] has been extended to the case of chaos control. Unlike the approach in [31], the analytical stability conditions for chaos control have been incorporated within the optimisation algorithm itself. This ensures stability of the optimised solutions in all cases, even when considering different initial conditions. To the best of our knowledge, the present paper can be considered as the first attempt for active control policy design for commensurate and incommensurate FO chaotic systems in a multi-objective framework with guaranteed analytical stability considerations.

The rest of the paper is organised as follows. Section 2 outlines the preliminary background of fractional calculus along with the numerical methods for simulating FO chaotic systems. Section 3 introduces the FO financial system and proposes the mathematical underpinnings of the active control strategy. Section 4 highlights the need for multi-objective optimisation in chaos control and describes the NSGA-II, ev-MOGA and MOEA/D algorithms briefly as multi-objective optimisers. Section 5 illustrates the results and discussions. The paper ends with the conclusions in Section 6 followed by the references.

## 2. Mathematical preliminaries
### 2.1. Basics of fractional calculus

Fractional calculus is an extension of the integer order successive differentiation and integration for any arbitrary real order. The fundamental operator representing the non-integer order differentiation or integration is given by ${}_aD_t^\alpha$ in (1), where $\alpha \in \mathbb{R}$ is the order of the differ-integration and $a$ and $t$ are the bounds of the operation.

$$_aD_t^\alpha = \begin{cases} d^\alpha/dt^\alpha, & \alpha > 0 \\ 1, & \alpha = 0 \\ \int_a^t (d\tau)^\alpha, & \alpha < 0 \end{cases} \tag{1}$$

There are three main definitions of fractional calculus *viz.* the Grünwald-Letnikov (GL), Riemann-Liouville (RL) and Caputo. Other definitions like that of Weyl, Fourier, Cauchy, Abel and Nishimoto also exist. In the FO systems and control related literatures mostly the Caputo's fractional differentiation formula is referred. This typical definition of fractional derivative is generally used to derive FO transfer function models from FO ordinary differential equations with zero initial conditions. According to Caputo's definition, the $\alpha^{th}$ order derivative of a function $f(t)$ with respect to time is given by (2)

$$_aD_t^\alpha f(t) = \frac{1}{\Gamma(m-\alpha)} \int_a^t \frac{D^m f(t)}{(t-\tau)^{\alpha+1-m}} d\tau, \alpha \in \mathbb{R}^+, m \in \mathbb{Z}^+, m-1 \leq \alpha < m \tag{2}$$

where, $\Gamma(\alpha) = \int_0^t e^{-t} t^{\alpha-1} dt$ is the Euler's Gamma function. This definition is used in the present paper for realizing the fractional integro-differential operators of the chaotic system. The Caputo definition of fractional derivative is advantageous for control related applications over the Riemann-Liouville definition, since it only needs initial conditions for integer order derivatives and not initial conditions of fractional derivatives. The Laplace transform of the Caputo fractional derivative is given by (3) [32].

$$\mathcal{L}\left[_0D_t^\alpha f(t)\right] = \int_0^\infty e^{-st} {}_0D_t^\alpha f(t) dt = s^\alpha F(s) - \sum_{k=0}^{m-1} s^{\alpha-k-1} f^m(0), m-1 \leq \alpha < m \tag{3}$$

For zero initial condition, the Laplace transform of the three definitions boils down to the same expression $s^\alpha F(s)$ which are extensively used in many modern control applications. Also, continuous or discrete time rational approximation techniques for the FO differ-integrator $s^{\pm\alpha}$ are often employed for simulation [33].

### 2.2. Numerical method for simulating fractional order chaotic systems

Chaotic coupled differential equations can be numerically simulated using the power series expansion method, Adams-Bashford-Moulton predictor corrector method [34], continued fraction expansion (CFE) method [35] etc. As has been shown in Petras [32], the chaotic FO differential equations in (4) can be written in the form of a set of integral equations as in (7) and band limited rational approximations can be used for realizing the fractional differentials. This method is adopted in the present paper for simulating the FO chaotic system.

For a set of coupled fractional differential equations of the form (4),

$$\begin{aligned}
{}_0D_t^{q_1}x(t) &= f(x(t), y(t), z(t)) \\
{}_0D_t^{q_2}y(t) &= g(x(t), y(t), z(t)) \\
{}_0D_t^{q_3}z(t) &= h(x(t), y(t), z(t))
\end{aligned} \quad (4)$$

considering the fact that the fractional differ-integrals are linear operators, i.e.

$$_aD_t^\alpha (\lambda f(t) + \mu g(t)) = \lambda\, _aD_t^\alpha f(t) + \mu\, _aD_t^\alpha g(t) \quad (5)$$

and the fact that the FO derivative commutes with the integer order derivative, i.e.

$$\frac{d^n}{dt^n}\left(_aD_t^\alpha f(t)\right) = \,_aD_t^\alpha\left(\frac{d^n f(t)}{dt^n}\right) = \,_aD_t^{\alpha+n} f(t) \quad (6)$$

Equation (4) can also be written in the form of integral equations as (7).

$$\begin{aligned}
x(t) &= {}_0D_t^{1-q_1}\left(\int_0^t [f(x(t), y(t), z(t))]\,dt\right) \\
y(t) &= {}_0D_t^{1-q_2}\left(\int_0^t [g(x(t), y(t), z(t))]\,dt\right) \\
z(t) &= {}_0D_t^{1-q_3}\left(\int_0^t [h(x(t), y(t), z(t))]\,dt\right)
\end{aligned} \quad (7)$$

With the implementation of this transformation a Matlab/Simulink based environment is capable of simulating the FO chaotic systems as shown by Petras [32].

Each value of the FO differ-integrals $\{1-q_1, 1-q_2, 1-q_3\}$ is rationalized with Oustaloup's 5[th] order rational approximation [36]. The FO differ-integrals are basically infinite dimensional linear filters. However, band-limited realisations of FO elements are necessary for simulation. In the present simulation study each FO element has been rationalised with Oustaloup's recursive filter [36] given by the equations (8) and (9). If it be assumed that the expected fitting range or frequency range of controller operation is $(\omega_b, \omega_h)$, then the higher order filter which approximates the FO element $s^\alpha$ can be written as (8) [37].

$$G_f(s) = s^\alpha = K\prod_{k=-N}^{N} \frac{s+\omega_k'}{s+\omega_k} \quad (8)$$

where the poles, zeros, and gain of the filter can be evaluated as:

$$\omega_k = \omega_b\left(\omega_h/\omega_b\right)^{\frac{k+N+\frac{1}{2}(1+\alpha)}{2N+1}},\ \omega_k' = \omega_b\left(\omega_h/\omega_b\right)^{\frac{k+N+\frac{1}{2}(1-\alpha)}{2N+1}},\ K = \omega_h^\alpha \quad (9)$$

In equations (8) and (9), $\alpha$ is the order of the differ-integration and $(2N+1)$ is the order of the filter. The present study considers a 5$^{th}$ order Oustaloup's rational approximation for the FO elements within the frequency range $\omega \in \{10^{-2}, 10^{2}\}$ rad/s [36], [37].

## 3. System description and theoretical formulation

The FO chaotic financial dynamical system [21] is given by (10).

$$\frac{d^{q_1} x}{dt} = z + (y-a)x$$
$$\frac{d^{q_2} y}{dt} = 1 - by - x^2 \qquad (10)$$
$$\frac{d^{q_3} z}{dt} = -x - cz$$

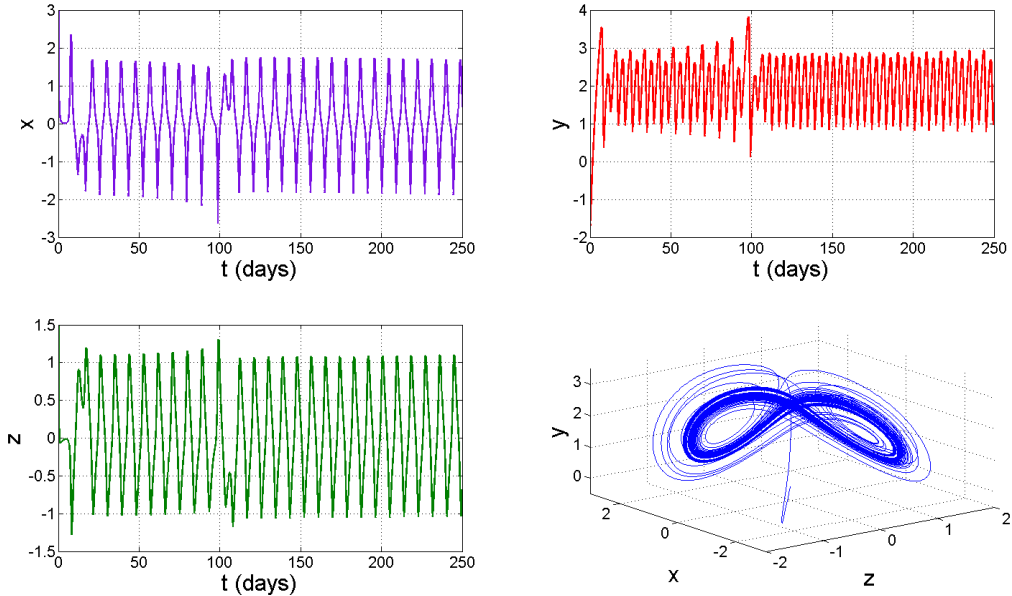

Figure 1: Phase portrait and state trajectories of the commensurate ($q_1$=$q_2$=$q_3$=q=0.9) FO financial system.

The state variables $x$, $y$, $z$ represent respectively the interest rate, the investment demand and the price index of a financial system. The first state variable ($x$) which is the interest rate can be influenced by the surplus between investment and savings along with structural adjustments from the prices. The second state variable ($y$) is in proportion to the rate of investment and an inversion with the cost of investment and interest rate. The third state variable ($z$) depends on the contradiction between supply and demand in the market and also gets influenced by the inflation rates [14]. The three constant coefficients $\{a, b, c\}$ represent the savings amount, the cost per investment and the elasticity of demand of commercial markets respectively. Figure 1 and Figure 2 show the variation in the three state variables i.e. interest rate, investment demand and price index with time (in days) [10] for the commensurate and incommensurate FO financial system respectively. All the three time series for both the systems exhibit erratic fluctuations, leading to a chaotic motion in the respective phase space diagrams. The second state variable (investment demand) shows more rapid fluctuation than the other two states indicating towards high spectral power in high frequency operation. Further details of the fractional financial system and its control are reported in Pan *et al.* [29]. Although stochastic modelling of financial systems exists in the literatures, for effective control policy development, the chaotic deterministic model of financial system is more popular.

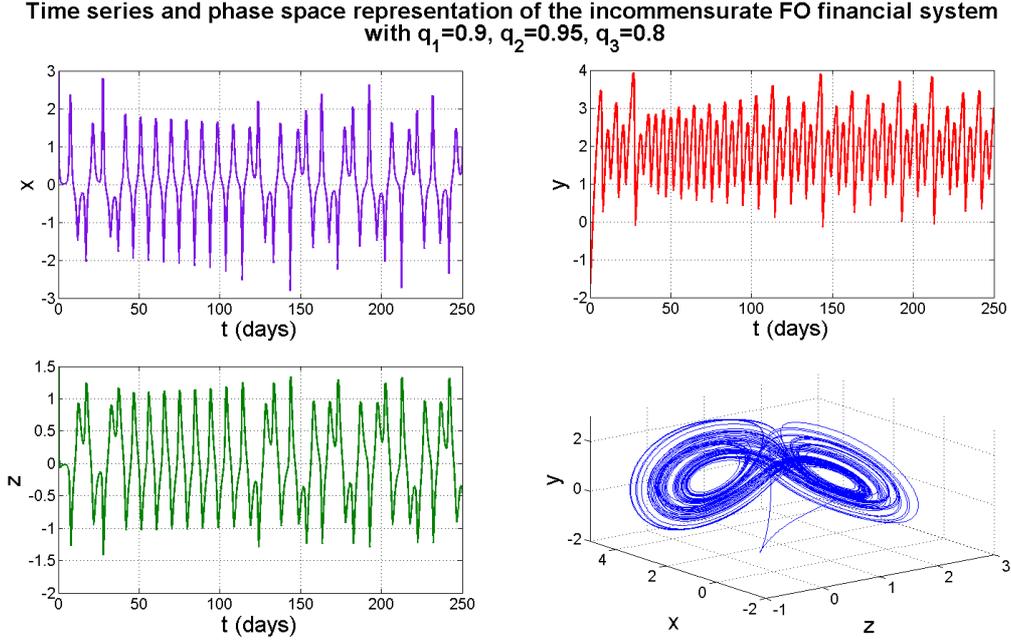

Figure 2: Phase portrait and state trajectories of the incommensurate ($q_1$=0.9, $q_2$=0.95, $q_3$=0.8) FO financial system.

For active control of the system described by (10), three active control functions $u_1(t), u_2(t), u_3(t)$ are considered to be applied in each of the three states of the system in (10) to yield the following set of equations.

$$\frac{d^{q_1} x}{dt} = z + (y-a)x + u_1(t)$$
$$\frac{d^{q_2} y}{dt} = 1 - by - x^2 + u_2(t) \quad (11)$$
$$\frac{d^{q_3} z}{dt} = -x - cz + u_3(t)$$

The nonlinear active state-feedback control functions are chosen as (12)-(14) in order to make the closed loop control system linear.

$$u_1(t) = V_1(t) - xy \quad (12)$$

$$u_2(t) = V_2(t) - 1 + x^2 \tag{13}$$

$$u_3(t) = V_3(t) \tag{14}$$

The terms $V_i(t) \, \forall i \in \{1,2,3\}$, are linear functions of the three system state variables $\{x, y, z\}$. Using equations (12), (13) and (14) in equation (11), we get (15).

$$\begin{aligned}
\frac{d^{q_1} x}{dt} &= z - ax + V_1(t) \\
\frac{d^{q_2} y}{dt} &= -by + V_2(t) \\
\frac{d^{q_3} z}{dt} &= -x - cz + V_3(t)
\end{aligned} \tag{15}$$

The active control terms $V_i(t) \, \forall i \in \{1,2,3\}$ can be represented by (16) where the constants $m_{ij} \in \mathbb{R}, \forall i, j \in \{1,2,3\}$.

$$\begin{bmatrix} V_1 \\ V_2 \\ V_3 \end{bmatrix} = \begin{bmatrix} m_{11} & m_{12} & m_{13} \\ m_{21} & m_{22} & m_{23} \\ m_{31} & m_{32} & m_{33} \end{bmatrix} \begin{bmatrix} x \\ y \\ z \end{bmatrix} \tag{16}$$

Thus (15) and (16) can be clubbed together to obtain (17).

$$D^q \begin{bmatrix} x \\ y \\ z \end{bmatrix} = P \begin{bmatrix} x \\ y \\ z \end{bmatrix} = \begin{bmatrix} -a + m_{11} & m_{12} & 1 + m_{13} \\ m_{21} & -b + m_{22} & m_{23} \\ -1 + m_{31} & m_{32} & -c + m_{33} \end{bmatrix} \begin{bmatrix} x \\ y \\ z \end{bmatrix} \tag{17}$$

where $q = [q_1, q_2, q_3]^T \in (0, 2)$.

The presence of squared and cross-product terms of the state variables in the active control inputs in (12)-(14) can be viewed as a nonlinear state feedback control design for the commensurate and incommensurate FO systems. The nonlinear control inputs makes the closed loop system linear to facilitate the established analytical stabilization schemes [32] for commensurate and incommensurate FO linear systems since such analytical stability design for the nonlinear counterpart is more involved and difficult to design.

The elements of the matrix in (16) are real. Hence it might be diagonalised to produce an equivalent control action where each of the control signals is a function of its own state only and not the other states. This would remove the number of couplings and reduce the complicacies in the design process. The number of elements would have been chosen as three instead of nine and this would reduce the burden of the optimisation algorithm. However, in many cases, the eigen-values obtained as a result of such diagonalisation might be complex conjugates and thus physical realisation of such state-feedback controllers with complex gains will be infeasible. Hence, all the nine components are chosen using optimization in (16) instead of choosing only the diagonals in the present approach. Another aspect of the control action that can be deduced from equations (12), (13), (14) is that, when the individual states become zero due to the application of the control action, i.e. $x = y = z = 0$, then $u_2$ becomes $(-1)$ and $u_1$ and $u_3$ becomes zero.

Now, in order to ensure asymptotic stability of system (17), the constants $m_{ij}$ must be chosen such that the eigen-values ($\lambda_k$) of matrix $P$ satisfy the following condition, known as Matignon's theorem [38].

$$\left|\arg\left(eig\left(P\right)\right)\right| = \left|\arg\left(\lambda_k\right)\right| > \frac{q\pi}{2}, 0 < q < 2 \tag{18}$$

for the commensurate FO system when $q_1 = q_2 = q_3 = q$.

For the incommensurate FO system, the asymptotic stability of Equation (17) can be derived using Deng's theorem as outlined in [39]. Let the incommensurate orders $q_i$ of equation (17) be written in the form $q_i = v_i/u_i$, $u_i, v_i \in \mathbb{Z}_+$. Let $m$ be the lowest common multiple (LCM) of $u_i$ and let $\gamma = 1/m$.

The characteristic equation can then be derived from (17) as shown in (19) by denoting the FO operators $s^{q_i}$ by $\lambda^{mq_i}$, where *diag* denotes a diagonal matrix [39], [40].

$$\det\left(s^{q_i}I - P\right) = \det\left(diag\left[\lambda^{mq_1} \quad \lambda^{mq_2} \quad \lambda^{mq_3}\right] - P\right) = \det\begin{bmatrix} \lambda^{mq_1} + a - m_{11} & -m_{12} & -1 - m_{13} \\ -m_{21} & \lambda^{mq_2} + b - m_{22} & -m_{23} \\ 1 - m_{31} & -m_{32} & \lambda^{mq_3} + c - m_{33} \end{bmatrix} = 0$$

(19)

If all the roots $\lambda_i, \forall i \in [1,2,3]$ of the characteristic equation as given in equation (19) satisfies $\left|\arg\left(\lambda_i\right)\right| > \frac{\gamma\pi}{2}, \forall i \in [1,2,3]$, then the system described by equation (17) is stable. The characteristic equation (19) is transformed to a higher integer-order polynomial equation if the incommensurate orders $q_i$ are considered as rational numbers. Matlab's Symbolic Math Toolbox function *solve()* has been used in the present paper to obtain the roots of the characteristic polynomial equation (19). Here, $m$ is 20 for the chosen incommensurate orders $q_1 = 0.9, q_2 = 0.95, q_3 = 0.8$. For checking the stability of incommensurate FO systems, it is sufficient to verify the argument of the roots in the primary Riemann sheet only, since the hyper-damped and ultra-damped roots in the higher Riemann sheets are always stable. Depending on the rational number representation of the incommensurate orders and by taking their LCM, the number of Riemann sheets could be tremendously high and also different roots might be distributed in different higher Riemann sheets which does not affect the stability of the incommensurate FO system. Therefore, in most literatures only the first Riemann sheet is considered for the stability checking as $(\pi/2m) \leq |\arg(\lambda_k)| \leq (\pi/m)$. Of course, one can place the root at specified locations, even in the higher Riemann sheets e.g. placement of ultra-damped roots as in the study by Bhalekar and Gejji [41], for the commensurate order systems. During the optimization based controller design for both the commensurate and incommensurate FO systems, a constraint is imposed such that at least two roots lie in the stable region of the primary Riemann sheet i.e. $(\pi/2m) \leq |\arg(\lambda_k)| \leq (\pi/m)$. This enables design of a relatively fast control system as opposed to the safe but very slow control operation when all the closed loop poles are pushed to the higher Riemann sheet (i.e. hyper-damped and ultra-damped roots) in the study by Das *et al.* [42].

## 4. Multi-objective optimisation for active control
### 4.1. Requirement for multi-objective optimisation

In any practical real world problem with constraints on resources, it is normal that satisfaction of one criterion to a greater extent would result in satisfaction of other conflicting criteria to a lesser extent. The concept of Pareto optimality [43] introduces this concept in economics on the topic of income distribution and economic efficiency. Consider that a finite number of goods are allocated among a set of individuals. If the economic allocation is a Pareto efficient one, then no individual can be made better off without one or more individuals being worse off. This same concept can be applied to the controller design problem for the financial system. There exists a trade-off in any controller design problem as illustrated in other literatures like [44], [45]. The two conflicting objective functions can be chosen as the Integral of the Time multiplied Squared Error (ITSE) ($J_1$) and the Integral of the Squared Deviation of Controller Output (ISDCO) ($J_2$). The two contradictory objectives can be mathematically expressed as

$$J_1 = ITSE_{set-point} = \int_{\psi}^{\infty} t e_{sp}^2(t) dt \qquad (20)$$

$$J_2 = ISDCO = \int_{\psi}^{\infty} (u(t) - u_{ss})^2 dt \qquad (21)$$

where, $\psi$ represents the instant of time, when the control signal was applied, $e_{sp}$ represents the error signal, i.e. the deviation of the chaotic trajectory from the desired set-point.

The first objective function $J_1$ tries to ensure fast tracking of the desired set-point. The time multiplication term assigns heavy penalty to the errors occurring at later stages and hence ensures faster settling time. The second objective function $J_2$ tries to reduce the change in the control signal as large control signal deviations necessitate large changes in the manipulated variables, which is not desirable [46]. This is because the manipulated variables are physical quantities and this might result in shocks to the system. $J_2$ is given by equation (21) and the term $\Delta u(t) = u(t) - u_{ss}$ represents the change in the absolute value of the control signal from its steady-state value. $J_1$ and $J_2$ are conflicting objectives since to reduce the steady state tracking error or to obtain fast tracking (i.e. to minimize $J_1$), the controller must exert more effort and hence the value of $J_2$ would increase and vice-versa. For evaluating equations (20) and (21), any of the three states $x, y, z$ or all of them taken together might be considered. The choice would depend on the practical constraints of the design. In this paper, the idea is to propose this concept of trade-off between different conflicting performance objectives and this is done considering only one of the states $y$ for evaluating equations (20) and (21). Therefore extensive simulations for all the possible cases are not reported in this paper.

In control design related optimization problems often conflict between two or more objectives are encountered i.e. speed of response and the control effort required [45]. But this is not obvious for any chosen control objective, for example several measures of tracking performance like Integral of Time multiplied Absolute Error (ITAE), ITSE, Integral of Absolute Error (IAE), Integral of Squared Error (ISE) may yield a Pareto front under a MOO framework but that does not guarantee that these objective are conflicting at all. The trade-off design becomes only important when there are physical constraints on arbitrarily increasing in one objective with an intention of keeping the other objectives at the same level. In Herreros *et al.* [44] several such conflicting objectives are formulated in frequency domain design of linear control systems. However, such frequency domain measures are not at all possible to derive for nonlinear chaotic systems. Therefore, for the present problem we had to rely on the time domain conflicting measures like ITSE and ISDCO [45]. From the point of view of the financial system this trade-off consideration is relevant as well. The three state variables $x$, $y$ and $z$ represent financial quantities which can be regulated, but cannot be arbitrarily increased or decreased in a practical setting. Depending on the way the objective functions are framed, one might want to supress the chaotic oscillations as quickly as possible and at the same time keep the deviation in these state variables and the required control effort to a minimum. This might be contradictory and hence a multi-objective methodology is required in the present scenario.

### *4.2. NSGA-II algorithm employed for multi-objective controller design*

A generalized multi-objective optimization framework can be defined as follows:

Minimize

$$F(x) = (f_1(x), f_2(x), ..., f_m(x)) \tag{22}$$

such that $x \in \Omega$.

where, $\Omega$ is the decision space, $\mathbb{R}^m$ is the objective space, and $F: \Omega \to \mathbb{R}^m$ consists of $m$ real valued objective functions.

Let, $u = \{u_1, ..., u_m\}$, $v = \{v_1, ..., v_m\} \in \mathbb{R}^m$ be two vectors, then $u$ is said to dominate $v$ if $u_i < v_i$ $\forall i \in \{1, 2, \cdots, m\}$ and $u \neq v$. A point $x^* \in \Omega$ is called Pareto optimal if $\nexists\ x\ |\ x \in \Omega$ such that $F(x)$ dominates $F(x^*)$. The set of all Pareto optimal points, denoted by PS is called the Pareto set. The set of all Pareto objective vectors, $PF = \{F(x) \in \mathbb{R}^m, x \in PS\}$, is called the Pareto Front. This implies that no other feasible objective vector exists which can improve one objective function without simultaneous worsening of some other objective function.

Multi-objective Evolutionary Algorithms (MOEAs) which use non-dominated sorting and sharing have higher computational complexity, uses a non-elitist approach and requires the specification of a sharing parameter. The NSGA-II removes these problems and is able to find a better spread of solutions and better convergence near the actual Pareto optimal front [47]. The pseudo code for the NSGA-II is as shown below [47], [48].

---

NSGA II Algorithm

Step 1: generate population $P_0$ randomly

Step 2: set $P_0 = (F_1, F_2, ...) = $ non-dominated-sort$(P_0)$

Step 3: for all $F_i \in P_0$

    crowding-distance-assignment$(F_i)$

Step 4: set t=0

    while (not completed)

        generate child population $Q_t$ from $P_t$

        set $R_t = P_t \cup Q_t$

        set $F = (F_1, F_2, ...) = $ non-dominated-sort$(R_t)$

        set $P_{t+1} = \phi$

        i=1

        while $|P_{t+1}| + |F_i| < N$

            crowding-distance-assignment$(F_i)$

            $P_{t+1} = P_{t+1} \cup F_i$

            i=i+1

        end

        sort $F_i$ on crowding distances

        set $P_{t+1} = P_{t+1} \cup F_i\left[1:(N - |P_{t+1}|)\right]$

        set $t = t + 1$

    end

return $F_1$

---

Here $N$ represents the number chromosomes in the population i.e. the population size. The NSGA-II algorithm converts $M$ different objectives into one fitness measure by composing distinct fronts which are sorted based on the principle of non-domination. In the process of fitness assignment, the solution set, not dominated by any other solutions in the population, is designated as the first front $F_1$ and the solutions are given the highest fitness value. These solutions are then excluded and the second non-dominated front from the remaining population $F_2$ is created and ascribed the second highest fitness. This method is iterated until all the solutions are assigned a fitness value. Crowding distances are the normalized distances between a solution vector and its closest neighbouring solution vectors in each of the fronts. All the constituent elements of the front are assigned crowding distances to be later used for niching. The selection is achieved in tournaments of size 2 according to the following logics.

a) If the solution vector lies on a lower front than its opponent, then it is selected.

b) If both the solution vectors are on the same front, then the solution with the highest crowding distance wins. This is done to retain the solution vectors in those regions of the front which are scarcely populated.

The population size is taken as 100 and the algorithm is run until the cumulative change in fitness function value is less than the function tolerance of $10^{-4}$ over 100 generations. The Crossover fraction is taken as 0.8 and an intermediate crossover scheme is adopted. The mutation fraction is taken as 0.2. For choosing the parent vectors based on their scaled fitness values, the algorithm uses a tournament selection method with a tournament size of 2. The Pareto front population fraction is taken as 0.7. This parameter indicates the fraction of population that the solver tries to limit on the Pareto front. The optimization variables are the components of the active control functions, i.e. $m_{ij} \ \forall i, j \in \{1,2,3\}$. Thus there are nine optimisation variables in total. To ensure that the solutions obtained are guaranteed to be stable, the stability criteria given by the Matignon's theorem for commensurate and incommensurate FO system are incorporated in the algorithm within each objective function evaluation. Thus the solutions that are generated through cross-over, mutation or reproduction in each generation, is first tested to see if they satisfy the stability criteria. In case the criteria are satisfied, the objective function is evaluated by simulating the chaotic system with the optimum controller gains, provided by the NSGA-II algorithm. In case the criteria are not satisfied then a high value of objective function is assigned to the solution without simulating the chaotic system since that particular controller cannot stabilize the system. This automatically assigns a fitness which is worse than the others to these unstable solutions. Therefore, over the generations, the algorithm rejects the unstable solutions and converges towards those regions in the solution space, which give stable controller values.

### 4.3. Testing other two MOO algorithms for the controller design: ev-MOGA and MOEA/D

For comparison with the obtained results of the NSGA-II algorithm, two other popular MOO algorithms *viz.* ev-MOGA and MOEA/D are also used for the active control policy design. The ε-MOGA is an elitist multi-objective evolutionary algorithm based on the concept of ε-dominance as discussed in Laumanns *et al.* [49]. In Deb *et al.* [50], a comparison of the ε-MOEA has been done with respect to other algorithms like NSGA-II, PESA (Pareto Envelope based Selection Algorithm), SPEA-II (Strength Pareto Evolutionary Algorithm-II) etc., where the ε-MOGA has been found to be superior. The ε-MOGA variable (ev-MOGA) is an improvisation over the ε-MOGA, which can characterise the Pareto front better than the ε-MOGA algorithm on several test-bench functions as indicated by Martínez-Iranzo *et al.* [51]. The ev-MOGA generates an ε-Pareto set ($\Theta_P^*$), which tries to converge towards the actual Pareto optimal set by adjusting the anchor points of the Pareto front $J(\Theta_P^*)$ dynamically and trying to preserve the ends of the Pareto front from being eliminated over the generations [52]. To achieve this, the objective space is divided into a fixed number of boxes, $n\_box_i$, specified by the user at the start of the algorithm. The algorithm consists of a main population of size $N_{ind_P}$, an auxiliary population of size $N_{ind_G}$ and an archive which stores the intermediate solutions and has an upper limit of $N_{ind\_max\_A} = 300$. The maximum number of generations is taken as $t_{max} = 2500$. In the present simulation $N_{ind_G}$ is taken as 10, $N_{ind_P}$ is taken as 30, $n\_box = \begin{bmatrix} 1000 & 1000 \end{bmatrix}$, the lower bounds of the variables is $\theta_{li} = \begin{bmatrix} -5 & -5 & -5 & -5 & -5 & -5 & -5 & -5 & -5 \end{bmatrix}$ and the upper bounds of the variables is $\theta_{ui} = \begin{bmatrix} 5 & 5 & 5 & 5 & 5 & 5 & 5 & 5 & 5 \end{bmatrix}$.

The MOEA/D algorithm as proposed in Zhang and Li [53], decomposes the MOO problem into a number of scalar optimization sub-problems and optimizes them simultaneously. The MOEA/D performs better than the NSGA-II algorithm on test bench problems and also has a lower computational complexity than the NSGA-II [53]. For the present exploration, the population size is taken as 70, the number of iterations as 1000 and the Tchebycheff approach is used for decomposition. The Tchebycheff approach associates a weight vector to each of the scalar sub-problems and different Pareto optimal solutions can be obtained by changing it. One of the pit falls of this approach is that its aggregation function is not smooth, but nevertheless, it can be used in the MOEA framework as there is no need to calculate the derivative of the aggregation function [53].

## 5. Simulations and results
### 5.1. Control of commensurate fractional order financial system

The system in (17) is simulated with the commensurate order $q_1 = q_2 = q_3 = q = 0.9$. The initial states $\{x_0, y_0, z_0\}$ of the system are chosen as $\{2, -1, 1\}$. Figure 3 shows the Pareto optimal set of solutions obtained from the NSGA-II algorithm. Each of the points represents a particular choice of the active control functions in (16) using three different MOO algorithm. The two axes denote the conflicting objectives of fast synchronisation performance and lower controller effort. It can be seen that each of these solutions satisfy the objectives to different extents. These solutions are non-dominated in the sense that it is not possible to find another set of solutions which would result in the performance improvement in both objectives. Thus there is a trade-off and a solution which would have a better performance in one of the objectives, would result in a lower performance in the other objective. Since the stability condition as given in (18) is checked as a sub-problem inside the MOO algorithms, these Pareto solutions are asymptotically stable. All the simulations reported here are based on a finite time window of $T = 50\,\text{sec}$.

**Table 1: Representative solutions on the Pareto front for the commensurate and incommensurate FO financial system**

| Class of FO model | Solution points | $J_1$ | $J_2$ | $m_{11}$ | $m_{12}$ | $m_{13}$ | $m_{21}$ | $m_{22}$ | $m_{23}$ | $m_{31}$ | $m_{32}$ | $m_{33}$ |
|---|---|---|---|---|---|---|---|---|---|---|---|---|
| Commensurate | $A_1$ | 1.005 | 11.627 | 2.300 | 0.607 | 1.633 | 1.357 | -2.133 | 1.041 | -0.663 | 0.818 | -0.129 |
| | $B_1$ | 9.921 | 1.248 | 2.049 | 0.902 | 0.708 | -2.412 | -0.747 | 0.000 | 0.572 | -0.033 | 0.757 |
| | $C_1$ | 97.152 | 1.125 | 2.044 | 0.948 | 0.583 | -2.449 | -0.543 | -0.049 | 0.649 | -0.038 | 0.897 |
| Incommensurate | $A_2$ | 1.001 | 15.995 | 4.999 | -0.018 | 3.312 | 3.259 | -4.880 | 0.357 | 0.959 | 4.934 | -0.920 |
| | $B_2$ | 2.224 | 1.477 | 3.124 | 0.188 | 2.251 | -0.625 | -0.445 | -1.673 | -0.517 | 0.214 | -1.250 |
| | $C_2$ | 38.298 | 1.053 | 2.975 | 0.485 | 1.838 | -1.089 | -0.255 | -1.925 | -1.803 | -0.413 | -1.385 |

**Table 2: Guaranteed stability for the three solutions on the Pareto front for two class of FO financial system**

| Class of FO model | Stability region (in degrees) | Solution points | Argument of the eigen-values of matrix P (in degrees) lying in the primary Riemann sheet | |
|---|---|---|---|---|
| Commensurate order q = 0.85 | $q\pi/2 = 81$ | $A_1$ | 81.0756 | -81.0756 |
| | | $B_1$ | 83.289 | -83.289 |
| | | $C_1$ | 81.1305 | -81.1305 |
| Incommensurate orders $q_1 = 0.9$, $q_2 = 0.95$, $q_3 = 0.8$ | $\pi/2m = 4.5$ | $A_2$ | -8.3496 | 8.3496 |
| | | $B_2$ | -7.1514 | 7.1514 |
| | | $C_2$ | 6.1084 | -6.1084 |

A set of representative solutions from the Pareto front is shown in Table 1. These solutions are the solutions at the extreme end of the Pareto front and the median solution. Table 1 reports the numerical values of the coefficients of the active control function in (16) along with the values of the two conflicting objectives as given in (20) and (21). In order to verify that the three chosen solutions (for both the commensurate and incommensurate order systems) are indeed obtained by satisfying the Matignon's stability criterion, the argument of the controlled system roots (at least two principal roots) are shown in

Table 2 to lie in the stable region of the primary Riemann sheet [33].

Figure 4 shows the time domain evolution of the states of the chaotic system along with the control input in the second state $u_2$. In all the cases, the time domain performance of the second state and the controller effort required in the second state is considered in the objective functions of equations (20) and (21), similar to that in [29]. As can be observed from Figure 4 that the *y*-state variable settles to the equilibrium point very quickly in the solution $A_1$. This is followed by solution $B_1$ and the solution $C_1$ is the most sluggish one amongst all the three and takes a long time to reach the equilibrium. In a practical setting, it is always desirable that the chaos is controlled in the shortest possible time. Thus while considering the time domain performance objective, the solution $A_1$ is the best and the solution $C_1$ is the worst. However, the opposite is observed for the controller effort. It can be seen that solution $A_1$ requires a large magnitude of control signal than solutions $B_1$ or $C_1$. In a real world case, the manipulated variable or the controller effort should be as small as possible as large changes in the manipulated variable might result in physical shock to the system or might not be possible to implement due to other constraints. Therefore solution $C_1$ is the best when considering the objective of lower control signal and the solution $A_1$ is the worst.

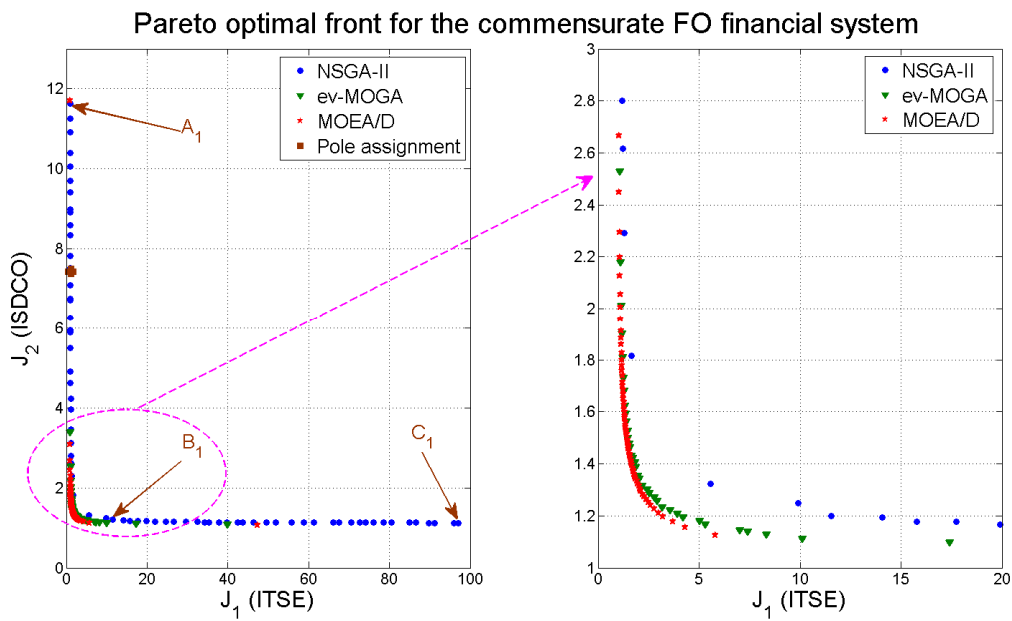

**Figure 3: Pareto optimal front for commensurate FO financial system.**

It is known that chaotic behaviour is observed in the FO financial system for specific values of the coefficients and the orders of the state variables ($q_1$, $q_2$, $q_3$) as reported in Chen *et al.* [21]. We have considered these specific values as the nominal financial system during the active nonlinear state feedback controller design since it exhibits chaos with the suggested parameters reported in [21]. However, in order to show that the proposed design technique works reliably for other fractional orders under same system structure and same co-efficients, the commensurate and incommensurate orders are decreased gradually. Our study shows that while gradually decreasing the fractional order, the chaotic behaviour disappears with $q = 0.8$ for commensurate FO financial system. It has been found from Figure 5 that the same controller parameters (median solution on the Pareto front) can suppress chaotic oscillations for other values of the commensurate fractional orders as well. They also satisfy the analytical stability conditions in this case. However, the controller parameters might not always be so robust for other systems as the issue of robustness has not been taken explicitly in the analytical formulation for the controller design. To design robust controllers, the same multi-objective methodology can be applied, but the mathematical stability formulation must be adequately changed.

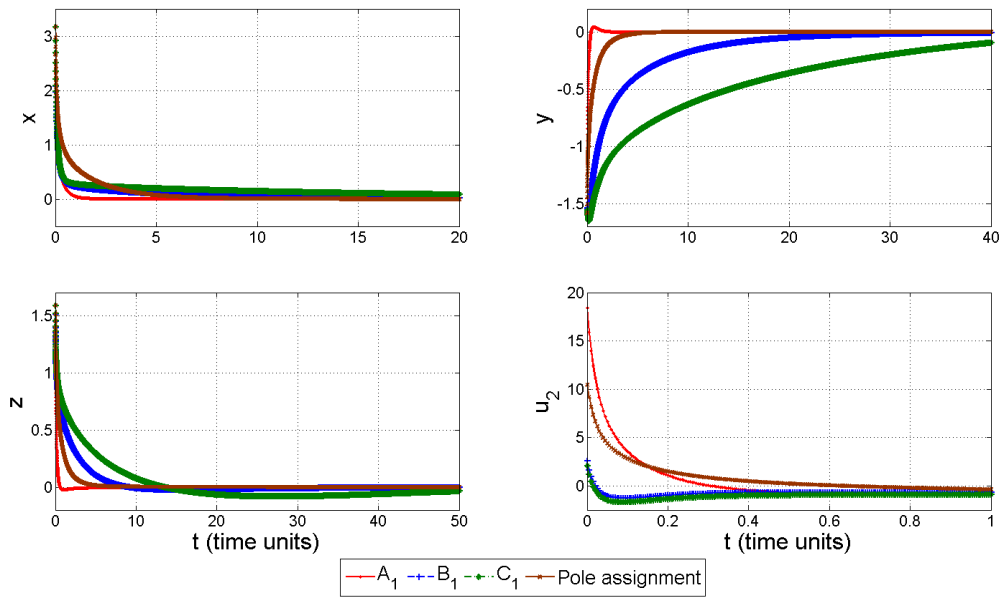

Figure 4: Chaos control with representative solutions on the Pareto front for the commensurate FO financial system

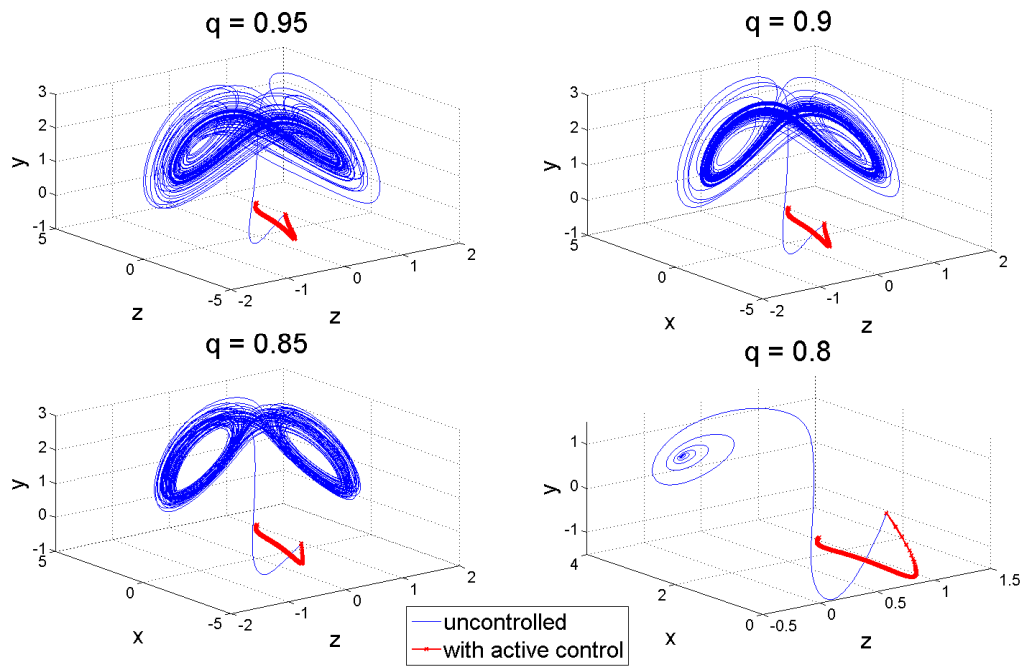

Figure 5: Uncontrolled and active controlled phase portraits for the change in commensurate FO financial system with the median solution ($B_1$) on the Pareto front.

### *5.2. Control of incommensurate fractional order financial system*

The system in (17) is simulated with the incommensurate order $q_1 = 0.9$, $q_2 = 0.95$, $q_3 = 0.8$. The initial states $\{x_0, y_0, z_0\}$ of the system are chosen as $\{2, -1, 1\}$. Figure 6 shows the Pareto optimal solutions for the incommensurate order case. The difference in simulation between this and the previous case is that the stability condition that is checked as a sub-problem in the MOO algorithms is different, as discussed before. Three representative solutions on the Pareto fronts (the ones at the extreme ends and the median solution) are reported in Table 1 and their corresponding time domain performances are shown in Figure 7. Among the two Pareto fonts in Figure 6 and Figure 3, it can be observed that the extent of the Pareto front is greater for the incommensurate order case. This implies that when the fractional orders of the chaotic system are different, then more variations in controller design are possible and the trade-offs that can be achieved among the different performance indices are greater.

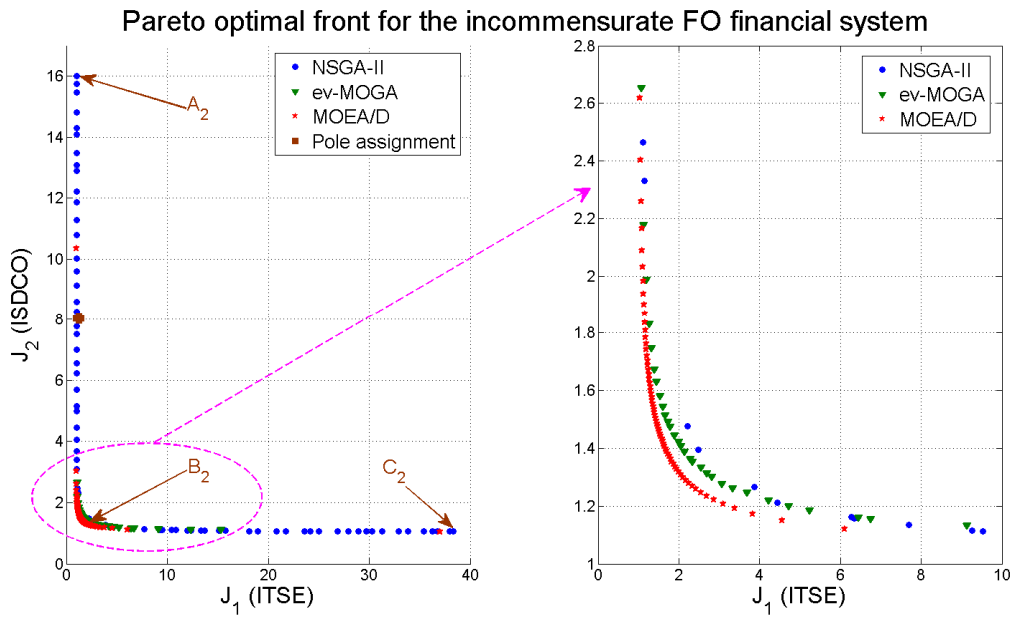

Figure 6: Pareto optimal front for incommensurate FO financial system.

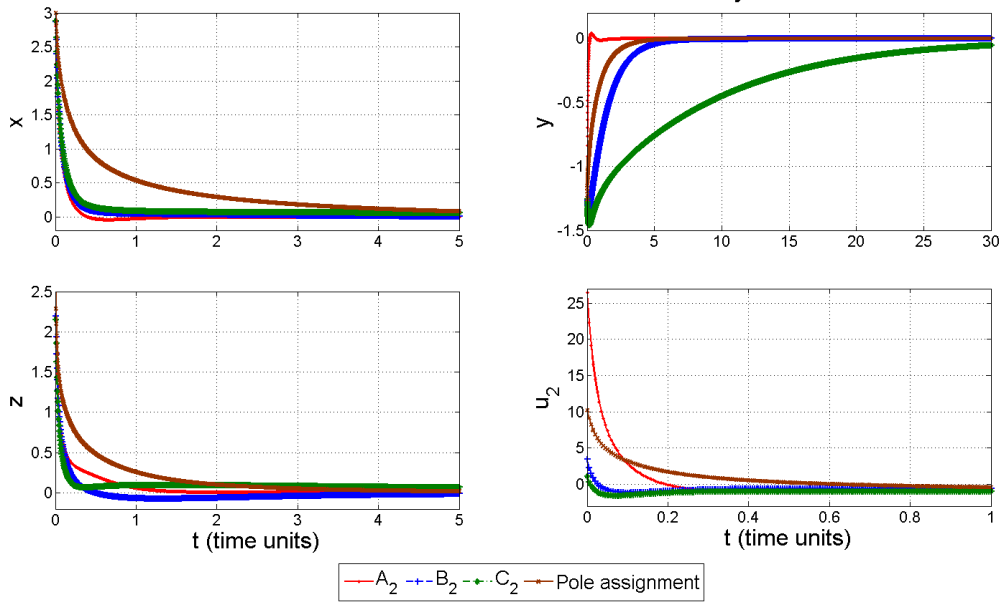

Figure 7: Chaos control with representative solutions on the Pareto front for the incommensurate FO financial system.

For both the commensurate and incommensurate FO financial systems, the Pareto fronts obtained by using three MOO algorithms *viz.* NSGA-II, ev-MOGA and MOEA/D are shown in Figure 3 and Figure 6 respectively. The simulations show that although the MOEA/D solutions dominate over that with the ev-MOGA and further NSGA-II (indicating that they are much fitter), the total spread of the Pareto fronts using either and ev-MOGA MOEA/D are much smaller. Therefore, in order to get a wide variety of solutions as a successful trade-off design one should consider the results with a larger Pareto front. We here report the time domain solution of the worst, best and median solutions of the Pareto front on either of the two control objectives. Similar exploration could be possible if the financial control policy designer gives more priority to the non-domination over the length of the Pareto front or the area covered by the Pareto front with respect to a chosen point in the dominated region [54], [55].

The nature of the solutions can be judged by Figure 6. It is further backed by the time domain evolution of the states of the chaotic system under the action of different controllers in Figure 7. The interpretations of the results are similar for the commensurate case. It can be observed from Figure 7, that the state *y* settles the fastest to the equilibrium point in solution $A_2$. This is followed by solutions $B_2$ and $C_2$ with the solution $C_2$ exhibiting the most sluggish response. However, the control effort required for the state *y*, is the least for solution $C_2$ and the maximum for solution $A_2$, with $B_2$ lying in between these two cases.

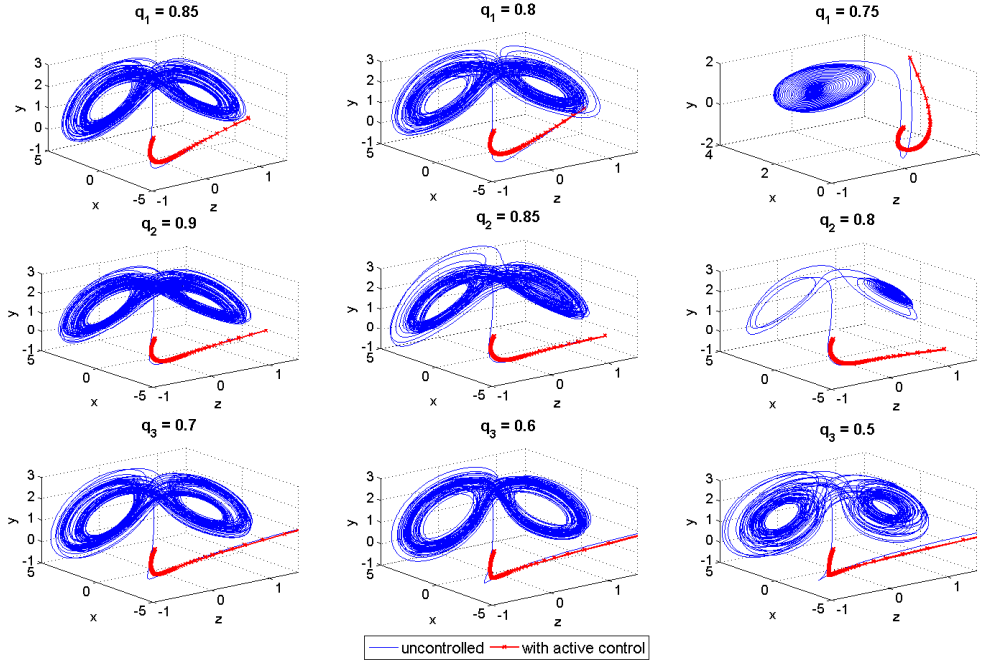

**Figure 8: Uncontrolled and active controlled phase portraits for gradual decrease in the incommensurate FO financial system orders from their nominal values q$_1$=0.9, q$_2$=0.95, q$_3$=0.8 with the median solution (B$_2$) on the Pareto front.**

For the incommensurate FO system, each of the three fractional orders of the three state variables are decreased gradually one at a time from their nominal values exhibiting chaos i.e. $q_1 = 0.9, q_2 = 0.95, q_3 = 0.8$, while keeping the other two fractional orders constant. In a scenario of separately changing the incommensurate orders the chaos disappears below $q_1 < 0.8, q_2 < 0.85, q_3 < 0.5$. However, for all the changes in the system order the same active control policy is capable of suppressing the chaotic oscillations for both the commensurate and incommensurate FO financial system which is evident from the phase portraits of the systems under active control in Figure 5 and Figure 8 respectively. In the present exploration, it is found that with the same controller, the system passes the analytical stability condition given by the Matignon's theorem. As expected the common control policy for the perturbed system would not be Pareto optimal, compared to different controllers tuned particularly for a fixed value of the commensurate and incommensurate orders. However these solutions still gives good results in terms of chaos suppression.

### 5.3. Comparison with direct pole assignment based active control approach

Most other controller design methods for FO chaotic systems, studied so far choose the active control functions heuristically such that the eigen-values are obtained in the stable region. Bhalekar and Gejji [41] have shown that an active control scheme could be designed intuitively such that all the eigen-values of the stability matrix becomes (-1) for a commensurate FO chaotic system. Following similar treatment for the commensurate order system, the $P$ matrix in (17) needs to be made a diagonal one such that the diagonal elements be the desired eigen-values. The corresponding choice of the constants $m_{ij}$ in (16) becomes

$$\left[ m_{ij} \right] = \begin{bmatrix} a+\bar{a} & 0 & -1 \\ 0 & b+\bar{b} & 0 \\ 1 & 0 & c+\bar{c} \end{bmatrix} \text{ in order to make the matrix } P = \begin{bmatrix} \bar{a} & 0 & 0 \\ 0 & \bar{b} & 0 \\ 0 & 0 & \bar{c} \end{bmatrix}.$$

Here the desired roots of the chaotic system under active control are chosen as $\bar{a} = \bar{b} = \bar{c} = -1$. This framework is applied on both the commensurate and incommensurate FO systems under study.

As expected such an intuitive choice will never produce optimum result as achieved by the proposed multi-objective optimization framework. In order to highlight this point, the obtained performance of the controlled system is shown in the two dimensional space between two conflicting objectives in Figure 3, using the method of active control design by Bhalekar and Gejji [41] for commensurate FO chaotic systems. Such an approach could not be extended easily for the incommensurate order systems since for such systems normally the stable roots within the primary Riemann sheet are only considered. Also depending on the incommensurate fractional orders and their LCM of associated denominators of those fractions $q_i$ there may be variable number of roots which is difficult to assign using a similar approach proposed for the commensurate order case in Bhalekar and Gejji [41].

Here, the LCM of $q_i$ is $m = 20$ for the chosen incommensurate orders of the system i.e. $q_1 = 0.9, q_2 = 0.95, q_3 = 0.8$. In order to get the characteristic equation in the form of desired roots, the choice of the constants in (16) would be $m_{11} = a+\bar{a}$, $m_{22} = b+\bar{b}$, $m_{33} = c+\bar{c}$, $m_{12} = m_{21} = m_{23} = m_{32} = 0$, $m_{13} = -1$, and $m_{31} = 1$ in order to make the stability

$$\text{matrix } P = \begin{bmatrix} -a+m_{11} & m_{12} & 1+m_{13} \\ m_{21} & -b+m_{22} & m_{23} \\ -1+m_{31} & m_{32} & -c+m_{33} \end{bmatrix} = \begin{bmatrix} -1 & 0 & 0 \\ 0 & -1 & 0 \\ 0 & 0 & -1 \end{bmatrix}.$$

With such a choice the desired pole locations results in $\bar{a}=\bar{b}=\bar{c}=-1$, the characteristic equation becomes as (23).

$$det\left[s^{q_i}I-P\right]=\det\left[\left(diag\begin{bmatrix}\lambda^{mq_1} & \lambda^{mq_2} & \lambda^{mq_3}\end{bmatrix}\right)-P\right]=\begin{bmatrix}\lambda^{18}+1 & 0 & 0\\ 0 & \lambda^{19}+1 & 0\\ 0 & 0 & \lambda^{16}+1\end{bmatrix} \quad (23)$$

As a result, the characteristic equation (23) yields several roots $\lambda_i$ which satisfies the incommensurate version of the Matignon's stability criterion. For the present choice of direct pole assignment scheme for the incommensurate FO system, all the roots of characteristic equation lies in the hyper-damped region i.e. $\left|\arg(\lambda_i)\right|>q\pi$ which signifies stable but slower operation of the system, compared to what can be achieved within a MOO framework. The reason for obtaining a dominated solution using the direct pole assignment technique compared to the Pareto solutions (in terms of control performance) is that the pole assignment scheme results in all hyper-damped and ultra-damped roots [33] leading to sluggish system operation and increase in the control effort. Unlike the commensurate FO system, in the incommensurate FO system stabilization problem, precise pole assignment becomes difficult due to the inherent higher order polynomial equation solving step. Although one can assign the eigen-values of the matrix *P* analytically, but manipulating the number and location of all the system roots for the characteristic equation which are distributed in different higher Riemann sheets for any arbitrary choice of incommensurate FO is still an open problem. Moreover, there would be only one solution by the direct pole-assignment approach for commensurate order system, as opposed to multiple solutions for the case of incommensurate FO system. As depicted in Figure 6, for incommensurate system also the direct pole assignment approach leads to a dominated solution, compared to what can be achieved by an MOO approach.

### 5.4. Discussions

It is also important to mention that FO controllers have been traditionally used for enhancing the robust stability properties of linear control systems. However, for nonlinear chaotic systems extension of the robust stability properties are expected to be more complex and have not been investigated yet, to the best of our knowledge, in the fractional calculus community itself. Moreover, the proposed control strategy does not use the concept of FO controller, rather it uses nonlinear state feedback control of FO systems. Although for both the commensurate and incommensurate order systems, the active control scheme designed under the nominal system parameters faithfully suppresses chaotic oscillations with gradual decrease in fractional orders and also passes the stability checking condition but during the controller design phase such a variation has not been considered. Therefore, although the same controller works well to stabilize different FO chaotic systems, this should not be confused with robust stability where the stability of all possible inter-mediate solutions are theoretically guaranteed and has been investigated for linear FO systems only in the past.

## 6. Conclusions

In this paper, an active control policy is derived for a FO chaotic financial system. The proposed method gives guaranteed stability, which is derived analytically for both the commensurate and incommensurate FO financial system. The active control functions are then chosen using three multi-objective evolutionary algorithms to satisfy two conflicting time domain performance objectives of fast settling to the equilibrium point and small amount of controller effort requirement. The comparison of three MOOs show that the NSGA-II yields the largest Pareto front over ev-MOGA and MOEA/D but a better non-dominated (although shorter) Pareto front could be achieved using MOEA/D. It is shown that the two design objectives cannot be simultaneously minimised using one particular controller. There exists a range of controllers on the Pareto front which satisfies one criterion better at the cost of performance deterioration in the other criterion. The designer can therefore choose a particular controller from these set of non-dominated solutions according to his specific problem requirements. The superiority of the proposed technique over the direct pole assignment approach [41] has also been highlighted. The effect of decreasing the fractional orders in the two type of systems (with the median solution of the controllers on the Pareto front) have been found to stabilize the chaotic systems and also pass the stability checking conditions Future work can be directed towards the multi-objective chaos control in the presence of uncertainty, noise etc. and extend the concept for robust stabilization scheme design for nonlinear chaotic systems.